\documentclass{article}
\usepackage{graphicx} 

\usepackage{amsthm,amssymb,amscd,amsmath,latexsym,indentfirst,color,amsfonts}
\usepackage[mathcal]{eucal}
\usepackage{geometry,amsthm,graphics,tabularx,shapepar,float}
\usepackage{stackrel}
\usepackage{extarrows}
\usepackage[all]{xypic}
\usepackage{tikz}
\usepackage{tikz-cd}
\usepackage{bm}
\usepackage{extarrows}
\usetikzlibrary{calc}
\theoremstyle{plain}
\theoremstyle{plain}
\newtheorem{theorem}{Theorem}[section]
\newtheorem{prop}[theorem]{Proposition}
\newtheorem{coro}[theorem]{Corollary}
\newtheorem{lemma}[theorem]{Lemma}
\newtheorem{ex}[theorem]{Example}

\newtheorem{rem}[theorem]{Remark}

\newtheorem{defi}[theorem]{Definition}

\newcommand{\ble}{\begin {lemma}}
\newcommand{\ele}{\end {lemma}}
\newcommand{\bde}{\begin {defi}}
\newcommand{\ede}{\end {defi}}
\newcommand{\bthm}{\begin {theorem}}
\newcommand{\ethm}{\end {theorem}}
\newcommand{\bco}{\begin {coro}}
\newcommand{\eco}{\end {coro}}
\newcommand{\bex}{\begin {ex}}
\newcommand{\eex}{\end {ex}}

\newcommand{\be}{\begin {equation}}
\newcommand{\ee}{\end {equation}}
\newcommand{\bp}{\begin {proof}}
\newcommand{\ep}{\end {proof}}
\newcommand{\bee}{\begin {equation*}}
\newcommand{\eee}{\end {equation*}}
\newcommand{\rt}{\rightarrow}
\newcommand{\lb}{\label}
\newcommand{\pt}{\partial}

\newcommand{\p}{\mathbf{p}}
\newcommand{\q}{\mathbf{q}}

\begin{document}
	\title{Normal forms of elements in the Weyl algebra and Dixmier Conjecture  }
	
\author{
	\small Gang Han \thanks{Corresponding author}\\
  \small	School of Mathematics, Zhejiang
\small	University\\
\small mathhgg@zju.edu.cn \\[1mm]	
\small		Zhennan Pan\\
\small School of Mathematics, Zhejiang
\small	University\\
\small pzn1025@163.com	\\[1mm]	
\small	Yulin Chen\\
\small School of Mathematics, Zhejiang
\small	University\\
\small tschenyl@163.com
}

\date{July 15, 2024}		
	\maketitle
	\begin{abstract}
		A result of A. Joseph says that any nilpotent or semisimple element $z$ in the Weyl algebra $A_1$ over some algebracally closed field $K$ of characterstic 0 has a normal form up to the action of the automorphism group of $A_1$. It is shown in this note that the normal form corresponds to some unique pair of integers $(k,n)$ with $k\ge n\ge 0$, and will be called the Joseph norm form of $z$. Similar results for  the symplectic Poisson algebra $S_1$ are obtained. The Dixmier conjecture can be reformulated as follows: For any nilpotent element $z\in A_1$ whose Joseph norm corresponds to 
		$(k,n)$ with $k>n\ge 1$, there exists no $w\in A_1$ with $ [z,w]=1$.  It is known to hold true  if $k$ and $n$ are coprime. In this note we show that the assertion also holds if $k$ or $n$ is prime. 
		Analogous results for  the Jacobian conjecture for $K[X,Y]$ are obtained.
	\end{abstract}
	
	

	\tableofcontents
	\section{Introduction}
	\setcounter{equation}{0}\setcounter{theorem}{0}
	
Let	$K$ be an algebraically closed field with characteristic 0.
Let	$A_1$ be the  Weyl Algebra generated over $K$ by $\p,\q $ with $[\q,\p]=1$. The polynomial algebra $K[X,Y]$ with the Poisson bracket defined by
\be \lb{}
\{f,g\}=\left |
\begin{array}{ll}
	{\pt f}/{\pt X} & {\pt f}/{\pt Y} \\
	{\pt g}/{\pt X} & {\pt g}/{\pt Y}
\end{array}
\right |,\ \ \  f,g\in K[X,Y], 
\ee
is a Poisson algebra. It is usually called the 1st symplectic Poisson algebra or the 1st canonical Poisson algebra, and will be denoted by $S_1$.
	
Let	$\text{Aut}(A_1)$ (resp. $\text{Aut}(K[X,Y])$,  $\text{Aut}(S_1)$) be the automorphism group of the $K$-algebra $A_1$ (resp.  the $K$-algebra $K[X,Y]$, the Poisson algebra $S_1$). 

	For $f,g\in K[X,Y]$, let $(f,g)$ denote the endomorphism $\phi$ of the $K$-algebra $K[X,Y]$ with $\phi(X)=f, \phi(Y)=g$ thus $\phi(F)=F(f,g)$ for $F\in K[X,Y]$. Let $J\phi=J(f,g)$ be the determinant  of the Jacobian matrix of $f$ and $g$, which is just $\{f,g\}$. For $\phi=(f,g)\in \text{Aut}(K[X,Y])$, $\phi\in \text{Aut}(S_1)$ iff  $J\phi=1$. 	It was shown by	L. Makar-limanov \cite{ml} that \be \lb{e11}\text{Aut}(A_1)\overset{}{\cong} \text{Aut}_n(K[X,Y]),\ee where $\text{Aut}_n(K[X,Y])$ is the group of automorphisms $(f,g) $ such that   $J(f,g)=1$. Thus  $\text{Aut}(A_1)$, $\text{Aut}(S_1)$ and $\text{Aut}_n(K[X,Y])$ are all isomorphic. 
	
	The Dixmier conjecture (abbr. $DC_1$) says, every endomorphism of the Weyl algebra $A_1$ is an automorphism. 
	Let $JC_2$ denote the Jacobian conjecture for $K[X,Y]$. It is known that $JC_2$ implies $DC_1$.
	
	In \cite{d}, the elements in $A_1$ which are not in the center are divided into  5 disjoint classes $\Delta_i, i=1,\cdots,5$, which are respectively strictly nilpotent, weak nilpotent, strictly semisimple, weak semisimple, and generic. This is usually called the Dixmier partition of $A_1$. The normal form of strictly nilpotent and strictly semisimple elements up to the action of $\text{Aut}(A_1)$ are well studied in \cite{d}. Then Joseph found a normal form for weak nilpotent and weak semisimple elements up to the action of $\text{Aut}(A_1)$, which are important for the current work. See Lemma 4.1 \cite{jo}. Let us introduce some notations before we state the first main result.
	
		For	$k,n\in \mathbb{Z}_{\geq 0}$, let $A_1(k,n)=\{z\in A_1|z = \sum a_{ij}\p^i\q^j,a_{kn} \neq 0$, and if $ i> k$ or $j> n$ then $a_{ij} = 0\}$. 	For $z,w\in A_1$, we write $z\simeq w$ if for some $ \phi \in \text{Aut}(A_1), w= \phi(z)$. 	
	 Let $G'_0$ be the subgroup of $\text{Aut}(A_1)$ consisting of those $\phi\in \text{Aut}(A_1)$ such that $\phi(\p)=a_1\p+c_1, \phi(\q)=a_2\q+c_2, a_1,a_2,c_1,c_2\in K,a_1a_2=1$. The following result is in Proposition \ref{p1}.
	 
	 	\begin{prop}
	 	Assume	$z\in A_1$  is nilpotent or semisimple, then there exists unique $(k,n)\in (\mathbb{Z}_{\geq 0})^2- \{ (0,0)\}, k\ge n\ge0,$ such that  $\phi(z)\in A_1(k,n)$ for some  $\phi \in \text{Aut}(A_1)$.  	If $z$ is strictly semisimple, then $(k,n)=(1,1)$.	 	
	 	If $z$ is nilpotent or weak semisimple, then  $k>n$; if $z\simeq z_1, z\simeq z_2$ and $z_1,z_2\in A_1(k,n)$, then $z_2=\phi(z_1)$ with  $\phi\in  G'_0$.
	 \end{prop}
	 The proof of this result uses a description of all the polynomial pairs $f$ and $g$ such that $(f,g)\in \text{Aut}(K[X,Y])$ by considering the Newton polygon of $f$ and $g$. See Proposition \ref{th1}.
For a	nilpotent or semisimple element $z\in A_1$, we will call the $z'\in A_1(k,n), k\ge n\ge 0,$ such that $z\simeq z'$   the Joseph normal form of $z$. In particular, the pair of integers $(k,n)$ with $k\ge n\ge 0$ is uniquely determined by $z$.

	 Analogous result are obtained for the Poisson algebra $S_1$. See Proposition \ref{p2}.
	
	Let  $\Gamma=\{z\in A_1|\exists w\in A_1, [z,w] = 1\}$.
	In Proposition \ref{p31}, we show that the Dixmier Conjecture $DC_1$ holds iff for any $k,n$ with $k>n\geq 1$, $ A_1(k,n)\cap \Gamma=\o$.
	 It is known that $ A_1(k,n)\cap \Gamma=\o$  if $k$ and $n$ are coprime. In this note we show the following result in Theorem \ref{main}.
	
		\bthm  Assume
	$ k> n \geq 1$. If $k$ or $n$ is prime, then $ \Gamma\cap A_1(k,n)=\o$.
	\ethm

	Analogous result for $JC_2$ are also obtained at the end of this note. \\

Below are some notations of this note.

For $f\in K[X,Y]-\{0\}$, assume $f=\sum a_{ij} X^iY^j$. Define $E(f)=\{(i,j)\in \mathbb{Z}^2| a_{ij}\ne 0\}$ to be the support of $f$.
For $(\rho,\sigma)\in \mathbb{R}^2- \{(0,0)\}$, we  set \[\mathbf{v}_{\rho,\sigma}(f)=sup\{\rho i+\sigma j|(i,j)\in E(f)\},\] which is usually called the $(\rho,\sigma)$-degree of $f$.
   Let
\[E(f;\rho,\sigma)=\{(i,j)\in E(f)|\rho i+\sigma j=\mathbf{v}_{\rho,\sigma}(f)\};\  \mathbf{f}_{\rho,\sigma}(f)=\sum_{(i,j)\in E(f;\rho,\sigma)} a_{ij} X^iY^j. \]

	A polynomial $f= \sum a_{ij}X^iY^j\in K[X,Y]$ such that  $a_{kn} \neq 0$  for some $k,n\ge0 $, and $a_{ij} = 0$ for $ i> k$ or $j> n$, is said to be of rectangular type $(k,n)$, or simply of  type $(k,n)$. A polynomial $f(X,Y)$ is said to be of rectangular type, if it is  of rectangular type $(k,n)$ for some $k,n\ge0 $.

	
	For $ f,g\in K[X,Y],$ we write $f\sim g$ if $ f = a g$ for some $ a \in K^{\times}$, where $K^{\times}=K-\{0\}$.

		For $0\ne f\in K[X,Y]$, $deg(f)$ denotes the total degree of $f$.

	Recall that the set of elements $\p^i \q^j,\ i,j\ge  0,$ is a basis for $A_1$.
	Let  $\phi\in \text{Aut}(A_1)$ with $\phi(\p)=f(\p,\q), \phi(\q)=g(\p,\q)$, where $f,g\in K[X,Y]$. Then we denote $\phi$ by 
	$(f(\p,\q),g(\p,\q))$.  
	
	One has the following $K$-linear isomorphism \be\lb{e00}  \Phi: A_1\rt K[X,Y], \p^i\q^j\mapsto X^iY^j,\ i,j\ge0.\ee
	
	For $z\in A_1-\{0\}$ and 
	$(\rho,\sigma)\in \mathbb{R}^2$ with $\rho+\sigma> 0$, let \[\mathbf{v}_{\rho,\sigma}(z)=\mathbf{v}_{\rho,\sigma}(\Phi(z)), \mathbf{f}_{\rho,\sigma}(z)=\mathbf{f}_{\rho,\sigma}(\Phi(z)).\]

	\section*{Acknowledgements}
	
	We would like to heartily thank Xiaoguang Wang for his help during this work. The proof of Proposition \ref{prop4} is due to him.
	
		\section{Normal forms of elements in $K[X,Y]$, $A_1$ and $S_1$  }
	\setcounter{equation}{0}\setcounter{theorem}{0}

First we recall the  Newton polygon of polynomials in $K[X,Y]$. The Newton polygon $NTP(f)$ of a polynomial  $f\in K[X,Y]-\{0\}$ is defined to be the convex hull in $\mathbb R^2$ of  $E(f)\cup \{O\}$, where $E(f)$ is the support of $f$ and $O$ is the origin of $\mathbb R^2$. For example, for $f=y+2x^3$, $NTP(f)$ is the triangle in 
$\mathbb R^2$  with vertices $(0,1), (3,0)$ and $O$.

For   $\phi=(f,g)\in \text{Aut}(K[X,Y])$, let 
\[ NTP(\phi)= (NTP(f),NTP(g)).\]
	
	Let $\tau =(Y,-X)\in \text{Aut}(K[X,Y])$. Then $\tau$ has order 4, $\tau^2=(-X,-Y) $ and $\tau^{-1} =(-Y,X)$. One has  For $\phi=(f,g)\in \text{Aut}(K[X,Y])$, 
	\[  \tau \circ (f,g)= (f(Y,-X), g(Y,-X)),\ \    (f,g) \circ\tau=(g,-f).\]

	It is clear then 	\[ NTP(\phi\circ\tau)= (NTP(g), NTP(f))\] and 	
	\[ NTP(\phi)= NTP(\phi\circ \tau^2)= NTP(\tau^2\circ\phi).\] Thus $ NTP(\tau^i\circ\phi\circ\tau^j), i,j\in \mathbb Z,$ equals one of $ NTP(\phi),  NTP(\tau\circ\phi),  NTP(\phi\circ\tau), NTP(\tau\circ\phi\circ\tau)$.

Let us recall two lemmas before we prove Proposition \ref{th1}, which describes  automorphisms of $K[X,Y]$ by their Newton polygons.
	
	\ble [Lemma 1.4 of \cite{d}]\lb{lm0}  Assume $(\rho,\sigma)\in \mathbb{R}^2- \{(0,0)\}$,  $f,g\in K[X,Y]- K$ are $(\rho,\sigma)$-homogeneous  and $\{f,g\}=0$. Assume $v=\mathbf{v}_{\rho,\sigma}(f), u=\mathbf{v}_{\rho,\sigma}(g) $ are both positive integers. Then 
	$f^u\sim g^v$. In particular if $s=u/v\in \mathbb Z$, then 	$f^s\sim g$.

	\ele 
	
	The following result is analogous to Lemma 2.7 of \cite{d}; see also Proposition 3.2 of \cite{jo}. For a proof see Proposition 2.2 of \cite{ht}.
	\ble\lb{lm5}
	
	Let $f,g\in K[X,Y]- \{0\}$. Assume that $(\rho,\sigma)\in \mathbb{R}^2- \{(0,0)\}$ with $\rho+\sigma> 0$. Set $f_*=\mathbf{f}_{\rho,\sigma}(f)$ and $g_*=\mathbf{f}_{\rho,\sigma}(g)$. Then
	
	(1) $\mathbf{f}_{\rho,\sigma}(fg)=f_*g_*$;
	
	(2) If $\{f_*,g_*\}\ne 0$ then \[\mathbf{f}_{\rho,\sigma}(\{f,g\})=\{f_*,g_*\}\ and \ \mathbf{v}_{\rho,\sigma}(\{f,g\})=\mathbf{v}_{\rho,\sigma}(f)+\mathbf{v}_{\rho,\sigma}(g)-(\rho+\sigma);\]
	
	(3) If $\{f_*,g_*\}= 0$, then $ \mathbf{v}_{\rho,\sigma}(\{f,g\})<\mathbf{v}_{\rho,\sigma}(f)+\mathbf{v}_{\rho,\sigma}(g)-(\rho+\sigma)$.
	
	\ele

	The following result is a corollary  of  Theorem 1 and 2 of \cite{di}, by considering the Newton polygon of  polynomials. 
	\begin{prop}\label{th1}
		Assume $\phi=(f,g)\in \text{Aut}(K[X,Y])$. 
		Then $\phi$  is of one of the following 4 types: 
		\begin{enumerate}
			\item 
			
			$\phi_1$ or $\tau\circ\phi_1$, where
			
			$\phi_1=(F_1,G_1):=(a_1X+a_0,b_2Y+b_0)$ , $a_i,b_i\in K,i=1,2,\ a_1,b_2\neq 0.$
			\item
			
			$\phi_2$ or  $\tau\circ\phi_2$ or $\phi_2\circ\tau$ or $\tau\circ\phi_2\circ\tau$, where
			
			$\phi_2=(F_2,G_2):=(a_1X+a_0,b_2Y+h(X)),\ \ a_1,b_2\neq0, h(X)\in K[X],\deg h=n\geq 1;$ in this case, $\mathbf{f}_{1,n}(G_2)=b_2Y+b_nX^n$, where $b_n$ is the leading coefficient of $h(X)$.

			\item $\phi_3=(F_3,G_3)$: $=(a_1X+a_2Y+a_0,b_1X+b_2Y+b_0), $
			$a_i,b_i\in K^\times, i=1,2, a_1b_2-a_2b_1\ne0.$
			\item 
			$\phi_4$ or  $\tau\circ\phi_4$ or $\phi_4\circ\tau$ or $\tau\circ\phi_4\circ\tau$, where
			
			$\phi_4=(F_4,G_4)$, such that for some $n\ge1$, $\mathbf{f}_{1,n}(F_4)=a(Y+cX^n)^k$,  
			$\mathbf{f}_{1,n}(G_4)=b (Y+cX^n)^{kl}, k,l\geq 1; a,b,c\neq 0.$ 
		\end{enumerate}
	\end{prop}
	\bp
	Assume $\phi=(f,g)\in \text{Aut}(K[X,Y])$.  If $deg(f)> deg(g)$ then we replace $(f,g)$ by $(g,-f)=(f,g)\circ\tau$, so one can assume that 
	$deg(f)\le deg(g)$. By Theorem 1 in \cite{di}, $NTP(f)$ is either a line segment (lying on the X-axis or Y-axis) or a triangle. 
	
	If $NTP(f)$ is  a line segment, then one has $f=a_1X+a_0, g=bY+h(X), a_1,b\ne 0$, or $(f,g)=\tau\circ (a_1X+a_0, bY+h(X))$. If $h(X)$ is a constant, then 
	$(f,g)$ is of the form 	$(F_1,G_1)$ or 	$\tau\circ(F_1,G_1)$; otherwise, 	$(f,g)$ is of the form 	$(F_2,G_2)$ or 	$\tau\circ(F_2,G_2)$.
	
	Now we assume that	$NTP(f)$ is a triangle. By  Theorem 1 in \cite{di}, one can assume that the vertices of $NTP(f)$ are $(m,0), (0,n), (0,0), m,n\ge1$, $m|n$ or $n|m$. In the case $n\nmid m$,  we replace $(f,g)$ by $\tau\circ(f,g)$. So we can assume  $n|m$. Let  $s=m/n$.  As $\{f,g\}=c\in K^\times$, $\{\mathbf{f}_{1,s}(f),\mathbf{f}_{1,s}(g)\}=c\ or\ 0$ by Lemma \ref{lm5}.
	
	
	(1)  If $\{\mathbf{f}_{1,s}(f),\mathbf{f}_{1,s}(g)\}\ne0$, let ${f}_*=\mathbf{f}_{1,s}(f), g_*=\mathbf{f}_{1,s}(g)$.Then $\{f_*, \{f_*,g_*\}\}=0 $. By Corollary 2.4 of \cite{jo}, there exists some  $\varphi\in \text{Aut}(K[X,Y])$, $\varphi(f_*)=dX^iY^j, d\ne0, i,j\ge0$. 
	
	If $s=1$ then one can choose $\varphi$ to be some  $(tX+uY,vX+wY), tw-uv=1$.
	If $s>1$ then one can choose $\varphi$ to be some $(X,Y+uX^s)$. 
	
	Assume  $\varphi(g_*)=h(X,Y)$, then $\{dX^iY^j,h\}=c\ne0$, thus $(i,j)=(0,1)\ or\ (1,0)$,  $\varphi(f_*)=dX\ or \ dY$.
	
	If $s=1$ then $f=a_1X+a_2Y+a_0, g=b_1X+b_2Y+b_0,  a_1b_2-a_2b_1\ne0$. As $NTP(f)$ is a triangle, $a_1,a_2\ne0.$ Thus $(f,g)$ is of the form 	$(F_3,G_3)$ if $b_1,b_2\ne0$, and  is of the form 	$(F_2,G_2)\circ\tau$ or 	$\tau\circ(F_2,G_2)\circ\tau$ if one of $b_1,b_2$ is 0.
	
	If  $s>1$ then as $\mathbf{f}_{1,s}(f)$ is not a monomial,  $\varphi(f_*)=dY$ and $f_*=d(Y-u	X^s)$. Then $g_*=bX, b\ne0$, and $(f,g)$ is of the form 	$(F_2,G_2)\circ\tau$. But now $deg(g)>deg(f)$, contradicts to the assumption 	$deg(f)\le deg(g)$.
	
	(2) Assume $\{\mathbf{f}_{1,s}(f),\mathbf{f}_{1,s}(g)\}=0$.	 
	
	By  Lemma \ref{lm0}, one has \be\label{eq4}f_*^{u}\sim g_*^{v}:\ f_*=\mathbf{f}_{1,s}(f), g_*=\mathbf{f}_{1,s}(g), v=\mathbf{v}_{n,m}(f), u=\mathbf{v}_{n,m}(g); v,u\in  \mathbb Z_{>0}.\ee Note that $(n,m)=n(1,s)$, $f_*^{u}=(f_*)^u $ and  $g_*^{v}=(g_*)^v$.
	
	As $deg(f)\le deg(g)$,  by  Theorem 2 in \cite{di} there exists $r_1>0, c_1\in K^*$, $deg(g-c_1f^{r_1})<deg(g).$ Then by (\ref{eq4}) one has \[g_*=c_1f_*^{r_1},  u/v=r_1\in \mathbb Z_{>0}.\]

	One has $(f,g-c_1f^{r_1})=(f,g)\circ(X,Y-c_1X^{r_1})\in \text{Aut}(K[X,Y])$. If $\{\mathbf{f}_{1,s}(f),\mathbf{f}_{1,s}(g-c_1f^{r_1})\}=0$ and $deg(f)\le deg(g-c_1f^{r_1})$, then one still has $\mathbf{f}_{1,s}(g-c_1f^{r_1})=c_2 f_*^{r_2} $ for some $c_2\ne0, r_2\ge1$, and continue this process until we find some $\tilde{g}=g-\sum_{i=1}^k c_if^{r_i}, r_i\ge1$, such that either $deg(f)>deg(\tilde{g})$ and  $\{\mathbf{f}_{1,s}(f),\mathbf{f}_{1,s}(\tilde{g})\}=0$, or $\{\mathbf{f}_{1,s}(f),\mathbf{f}_{1,s}(\tilde{g})\}\ne0$. Note that $(f,\tilde{g})\in \text{Aut}(K[X,Y])$.

	
	If $\{\mathbf{f}_{1,s}(f),\mathbf{f}_{1,s}(\tilde{g}) \}\ne0$, set $(f_1,g_1)=(f,\tilde{g})$ and stop the process. Set $t_0=1$, one has $f=f_1$ and $\mathbf{f}_{1,s}(f)\sim \mathbf{f}_{1,s}(f_1)^{t_0}$.
	
	If $\{\mathbf{f}_{1,s}(f),\mathbf{f}_{1,s}(\tilde{g}) \}=0$ and $deg(f)>deg(\tilde{g})$, let $(f_1,g_1)=(\tilde{g},-f)(=(f,\tilde{g})\circ\tau)$, and by the same reason as above, one has for some $t_0>0$, \[\mathbf{f}_{1,s}(f)\sim \mathbf{f}_{1,s}(f_1)^{t_0}.\] Now we continue the process in (2) for $(f_1,g_1)$. We will get $  (f_1,g_1),(f_2,g_2),\cdots,(f_l,g_l)$ in $\text{Aut}(K[X,Y])$ and $t_0,t_1,\cdots,t_{l-1}>0$  such that 
	\[\{\mathbf{f}_{1,s}(f_l),\mathbf{f}_{1,s}(g_l)\}\ne0;\  \mathbf{f}_{1,s}(f_i)\sim \mathbf{f}_{1,s}(f_{i+1})^{t_i}, i=0,1,\cdots,l-1, (f_0,g_0)=(f,g). \] Thus $\mathbf{f}_{1,s}(f)\sim \mathbf{f}_{1,s}(f_l)^{t}, t=t_0t_1\cdots t_{l-1}. $

	If $s>1$ then $(f_l,g_l)$ is of the form $(F_2,G_2)\circ\tau=(b_2Y+h(X),-a_1X-a_0),\ \ a_1,b_2\neq0, h(X)\in K[X],\deg h=s$. In this case, $(f,g)$ is of the form $(F_4,G_4)$ with \[\mathbf{f}_{1,s}(f)=a(Y+cX^s)^t,  
	\mathbf{f}_{1,s}(g)=b(Y+cX^s)^{tr_1}, t\geq 1, r_1=u/v,a,b,c\neq 0, b=a^{r_1}c_1.\]
	
	If $s=1$ then $(f_l,g_l)$ is of the form $(F_3,G_3)=(a_1X+a_2Y+a_0,b_1X+b_2Y+b_0), a_1,a_2\ne0$.
	In this case, $(f,g)$ is of the form $(F_4,G_4)$ with \[\mathbf{f}_{1,1}(f)=(a_1X+a_2Y)^t,  
	\mathbf{f}_{1,1}(g)=c_1\cdot(a_1X+a_2Y)^{tr_1}, t\geq 1,a_1,a_2,c_1\neq 0.\]
	
	\ep
	The above result can also be proven by induction using the decomposition of $\text{Aut}(K[X,Y])$ as the 
	amalgamated product of the group of affine automorphisms  of $K[X,Y]$ and the group of triangular automorphisms of $K[X,Y]$ over their intersections.
	
	Here is a corollary of above result.
	
	\bco\lb{c1}
	Assume $f\in K[X,Y]$ such that $(f,g)\in \text{Aut}(K[X,Y])$ for some $g\in K[X,Y]$. Then there are 2 possibilities.
	
	(1) $NTP(f)$ is a line segment, and $f$ or $\tau(f)$ is $aX+b, a\ne0$;
	
	(2) $NTP(f)$ is a triangle, and for some $n\ge1$, $\mathbf{f}_{1,n}(f)$ or $\mathbf{f}_{1,n}(\tau(f))$ equals $a(Y+cX^n)^k, k\geq 1; a,c\neq 0.$ 
	\eco
	
Let	$G_0$ denote the subgroup of $ \text{Aut}(K[X,Y])$ consisting of all the automorphisms\\ $(a_1X+c_1,a_2Y+c_2)$ with  $ a_i,c_i\in K,i=1,2;\ a_1,a_2\neq 0$. 
	
	\ble\label{lm2}
	Assume $h_i\in K[X,Y]$ and $h_i$ is of type  $(k_i,n_i)$, $ k_i>n_i\ge0,\  i=1,2.$ Assume  there exists $ \phi\in \text{Aut}(K[X,Y])$ such that $h_2=\phi(h_1).$ Then $$(k_2,n_2)=(k_1,n_1),$$ and $\phi$ can be chosen in the group  $G_0.$
	
	\ele
	\begin{proof}
		Assume first $k_1> n_1\geq 1$. Assume
		$h_1(X,Y)= \sum a_{ij}X^iY^j\in K[X,Y]$, $a_{k_1n_1} \neq 0$, and $a_{ij} = 0$ for $ i> k_1$ or $j> n_1$. Then $h_1(X,Y)=\lambda_1X^{k_1}Y^{n_1}+\widetilde{h_1}(X,Y)$ with $\lambda_1=a_{k_1n_1}$.	 The proof below uses the notations in  Proposition \ref{th1}.

		If $\phi= (F_4,G_4) $ defined in (4) of Proposition \ref{th1}, then $ \phi(h_1)=\lambda_1 F_4^{k_1}G_4^{n_1}+\widetilde{h_1}(F_4,G_4) $ and it is clear that  $$\mathbf{f}_{1,n}(\phi(h_1))=\lambda_1\mathbf{f}_{1,n}(F_4^{k_1}G_4^{n_1})=\lambda_1 a^{k_1} b^{n_1}(Y+cX^n)^{r}, r=kk_1+kln_1>0; a,b,c\ne0.$$ Thus $h_2=\phi(h_1)$ is not of rectangular  type,  contradicts to the assumption. Similarly  $\phi$ cannot be $\tau\circ \phi_4,\phi_4\circ\tau$ or $\tau\circ\phi_4\circ\tau$ in (4) of Proposition \ref{th1}.
		
		Similarly if $\phi=(F_3,G_3)$ in (3) of Proposition \ref{th1}, then $\mathbf{f}_{1,1}(\phi(h_1))$ is not a monomial,  contradicts to the assumption.
		
		If $\phi = (F_2,G_2)$ in (2) of Proposition \ref{th1}, and assume $b_nX^n$ is the leading term of $h(X),\ n\geq 1$, then we have 
		$$\mathbf{f}_{1,n}(\phi(h_1))=\lambda_1(a_1X)^{k_1}(b_2Y+b_nX^n)^{n_1}.$$
		As $n_1\geq 1$, $h_2$ is not of rectangular type, contradicts to the assumption.
		Similarly, $\phi$ can not be $\tau\circ \phi_2,\phi_2\circ\tau$ or $\tau\circ\phi_2\circ\tau$ in (2) of Proposition \ref{th1}.
		
		If $\phi = (F_1,G_1)$, then $h_2=\phi(h_1)$ is of rectangular type  $(k_1,n_1)$, thus		
		$(k_2,n_2)=(k_1,n_1)$, and $\phi$ is in $G_0$.	If $\phi =\tau\circ (F_1,G_1)$, then $h_2=\phi(h_1)$ is of rectangular type  $(k_2,n_2)=(n_1,k_1)$, contradicts to $k_2> n_2$.	
		
		
		Next assume $k_1> n_1=0.$ Now $h_1$ is a polynomial in $X$. Reasoning as above we find that $\phi(h_1)$ satisfies the requirement iff $\phi=\phi_1=(F_1,G_1)$ or $\phi=\phi_2=(F_2,G_2)$ in  Proposition \ref{th1}. Now one still has $(k_2,n_2)=(k_1,n_1)$. As $\phi_2(h_1)=\phi_1(h_1) $, thus $\phi$ can be chosen in $G_0$.
		
	\end{proof}
	
	\bco\label{c3}
	Assume $h_i\in K[X,Y]$ and $h_i$ is of  type $(k_i,n_i)\ne (0,0)$, $ k_i,n_i\ge0 ,i=1,2.$ Assume  there exists $ \phi\in \text{Aut}(K[X,Y])$ such that $h_2=\phi(h_1).$ Then $$(k_2,n_2)=(k_1,n_1)\ or\ (k_2,n_2)=(n_1,k_1),$$ and $\phi$ can be chosen in the group
	$C_4 G_0$, where $C_4$ is the order 4 cyclic group generated by  $\tau=(Y,-X)$ and  $G_0$ is normalized by $C_4$.
	\eco
	\bp
	Suppose $k_1\neq n_1$. Then one can assume $k_1>n_1$. As $h_1$ and  $h_2=\phi(h_1)$ are both of rectangular type,  by lemma \ref{lm2}, one can require $\phi$ to  be of the form $\phi_1$ or $ \tau\circ\phi_1$, where $\phi_1\in G_0.$ If $\phi=\phi_1$, then $(k_2,n_2)=(k_1,n_1)$; if $\phi= \tau\circ\phi_1$, then $(k_2,n_2)=(n_1,k_1)$.
	
	
	Suppose $k_1=n_1\geq 1$,by the same reasoning as in Lemma \ref{lm2} one finds that $\phi$ must also be of the form $\phi_1$ or $\tau\circ\phi_1$ with  $\phi_1\in G_0$, and one always has $(k_2,n_2)=(k_1,n_1)$.
	\ep

		One knows that $\text{Aut}_n(K[X,Y])=L\underset{H}{*}T$ is the amalgamated product of $L$ and $T$ over $H=L\cap T$, where
	$$L=\left\{(a_{11}X+a_{12}Y+b_{1},a_{21}X+a_{22}Y+b_{2})\bigg|a_{11}a_{22}-a_{12}a_{21}=1\right\}$$
	is the subgroup of affine automorphisms of $A_1$, 
	$$T=\left\{(aX+b,a^{-1}Y+f(X))|f\in K[X],a\in K^\times,b\in K \right\}$$ is the subgroup of triangular automorphisms of $K[X,Y]$, and 
	$$H=L\cap T=\left\{(aX+b,a^{-1}Y+c_{1}X+c_{2})|c_{1},c_{2},b\in K,a\in K^\times \right\}.$$

	Let
	$$S_{1}=\left\{(X+aY,Y)\big|a\in K\}\cup\{(Y,-X)\right\}$$
	be a set of left coset representatives of $L/H$. let
	$$S_{2}=\left\{(X,Y+f(X))\big|f\in K[X],\deg f\geq2,or\ f=0\right\}$$
	be a set of left coset representatives of $T/H$. Then any $\psi\in \text{Aut}_n(K[X,Y])$ can be uniquely written as
	\be\lb{e12}\psi=\varphi_{1}\varphi_{2}\cdots\varphi_{n}\xi: \ \xi\in H, n\ge0, \varphi_i\in S_{t_i}-\{1\}, t_i\in\{1,2\}\ and\ t_i\ne t_{i+1}, i=1,\cdots, n-1. \ee
	
	One has the embedding \[\text{Aut}(K[X,Y])\rt K[X,Y]\times K[X,Y], \phi=(f,g)\mapsto (f,g) \] and \[\text{Aut}(A_1)\rt A_1\times A_1, \phi=(f(\p,\q),g(\p,\q))\mapsto (f(\p,\q),g(\p,\q)). \]
	The map \[K[X,Y]\times K[X,Y]\rt A_1\times A_1, (f,g)\mapsto (f(\p,\q),g(\p,\q)) \] thus induces $\eta:L\cup T\rt \text{Aut}(A_1), \phi=(f,g)\mapsto \phi'=(f(\p,\q),g(\p,\q))$, since it is clear that $\eta(\phi)\in \text{Aut}(A_1)$ for $\phi\in L\cup T$.  Let $L'=\eta(L), T'=\eta(T), H'=\eta(H), S'_{i}=\eta(S_i), i=1,2 $.
	
	By (\ref{e11}) one has
	 $\text{Aut}(A_1)=L'\underset{H'}{*}T'$, and the element $\psi$ in (\ref{e12}) corresponds to  	$\psi'=\varphi'_{1}\varphi'_{2}\cdots\varphi'_{n}\xi'$, 	with $\xi'=\eta(\xi)\in H',\varphi'_{i}=\eta(\varphi_i)\in S'_{t_i}-\{1\}$. So
each element $\kappa$ in $\text{Aut}(A_1)$ can be written as $\psi_1\delta_1\psi_2\delta_2\cdots \psi_k\delta_k \xi, or\  \psi_1\delta_1\psi_2\delta_2\cdots \psi_k\delta_k \psi_{k+1}\xi,$ $k\ge0, \psi_i\in S_1'-\{1\}, \delta_i\in S_2'-\{1\}, \xi\in H'$ and we allow $\psi_1$ to be 1. As $\xi, \psi_{k+1}\xi\in L'$, $\kappa= \psi_1\delta_1\cdots \psi_k\delta_k \zeta, k\ge0, \zeta\in L'$.  Then by induction on $k$, one can prove  the following result analogous to Proposition \ref{th1}.
	 	 
	 	\bco\label{}
	 Assume $\phi=(f(\p,\q),g(\p,\q))\in \text{Aut}(A_1)$ with  $f,g\in K[X,Y]$. Then $(f,g)$  is of one of the following 4 types: 
	 \begin{enumerate}
	 	\item  	$\phi_1$ or $\tau\circ\phi_1$, where
	 	
	 	$\phi_1=(F_1,G_1):=(a_1X+c_1,a_2Y+c_2)$ , $a_i,c_i\in K,i=1,2,\ a_1a_2=1.$
	 	\item  	$\phi_2$ or  $\tau\circ\phi_2$ or $\phi_2\circ\tau$ or $\tau\circ\phi_2\circ\tau$, where		
	 	
	 	$\phi_2=(F_2,G_2):=(a_1X+c_1,b_2Y+h(X)),\ \ a_1b_2=1, h\in K[X],\deg h=n\geq 1;$ in this case, $\mathbf{f}_{1,n}(G_2)=b_2Y+b_nX^n$.
	 	\item 		
	 	$\phi_3=(F_3,G_3)$: $=(a_1X+b_1Y+c_1, a_2X+b_2Y+c_2), $
	 	$a_i,b_i\in K^\times, i=1,2,\  a_1b_2-a_2b_1=1.$
	 	\item 	$\phi_4$ or  $\tau\circ\phi_4$ or $\phi_4\circ\tau$ or $\tau\circ\phi_4\circ\tau$, where
	 	
	 	$\phi_4=(F_4,G_4)$		 such that for some $n\ge1$, $\mathbf{f}_{1,n}(F_4)=a(Y+cX^n)^k$,  
	 	$\mathbf{f}_{1,n}(G_4)=b (Y+cX^n)^{kl}, k,l\geq 1; a,b,c\neq 0.$ 
	 	
	 \end{enumerate}
	 \eco

	 One also has the following result analogous to Corollary \ref{c1}.
	  \bco
	 Assume $f\in K[X,Y]$ such that $(f(\p,\q),g(\p,\q))\in \text{Aut}(A_1)$ for some $g\in K[X,Y]$. Then there are 2 possibilities.
	 
	 (1) $NTP(f)$ is a line segment, and, $f$ or $\tau(f)$ equals $aX+b, a\ne0$;
	 
	 (2) $NTP(f)$ is a triangle, and for some $n\ge1$, $\mathbf{f}_{1,n}(f)$ or $\mathbf{f}_{1,n}(\tau(f))$ equals $a(Y+cX^n)^k, k\geq 1; a,c\neq 0.$ 
	 
	 \eco

		Recall that for	$k,n\in \mathbb{Z}_{\geq 0}$, $A_1(k,n)=\{z\in A_1|z = \sum a_{ij}\p^i\q^j,a_{kn} \neq 0$, and if $ i> k$ or $j> n$ then $a_{ij} = 0\}$.
	
	Let $\tau'=(q,-p)$ and $G'_0=\{(a_1\p+c_1,a_2\q+c_2)|a_1,a_2,c_1,c_2\in K,a_1a_2=1\}$, which is a subgroup of $\text{Aut}(A_1)$. 
	Analogous to Lemma \ref{lm2} and Corollary \ref{c3}, one has
	\begin{prop}\label{p5}
		Assume $z_i\in A_1(k_i,n_i),k_i,n_i\geq 0,i=1,2.$ Assume there exists $ \phi\in \text{Aut}(A_1) $ such that $ z_2=\phi(z_1). $ Then 
		
		(1) 	If  $k_i>n_i\ge 0$ for $i=1,2$, then $(k_2,n_2)=(k_1,n_1)$ and $\phi$ can  be chosen in the group $G'_0$.
		
		(2) In general, one has 
		$$(k_2,n_2)=(k_1,n_1)\ or \ (n_1,k_1),$$ and
		$\phi$ can  be chosen in the group $C'_4\cdot G'_0$, where $C'_4$ is the order 4 cyclic group generated by the automorphism $\tau'$ and $G'_0$ is normalized by $C'_4$.
		
	\end{prop}

	\begin{theorem}[Lemma 4.1 of \cite{jo}]\label{th6}
		Let	$z\in A_1.$ Assume $z$ is nilpotent or semisimple, then there exists $\phi \in \text{Aut}(A_1)$ such that $\phi(z)\in A_1(k,n)$, where $k,n\in \mathbb{Z}_{\geq 0}, (k,n)\neq (0,0).$
			If $z$ is strictly semisimple, then $(k,n)=(1,1)$; otherwise $k\ne n$, and one can choose $k,n$ such that $k>n$.
	\end{theorem}
	
		For $z,w\in A_1$, we write $z\simeq w$ if for some $ \phi \in \text{Aut}(A_1), w= \phi(z)$.

	Combine Theorem \ref{th6} and Proposition \ref{p5}, one has
	\begin{prop}\lb{p1}
		Assume	$z\in A_1$  is nilpotent or semisimple, then there exists unique $(k,n)\in (\mathbb{Z}_{\geq 0})^2- \{ (0,0)\}, k\ge n\ge0,$ such that  $\phi(z)\in A_1(k,n)$ for some  $\phi \in \text{Aut}(A_1)$. 		
		If $z$ is strictly semisimple, then $(k,n)=(1,1)$.		
		If $z$ is nilpotent or weak semisimple, then  $k>n$; if in addition $z\simeq z_1, z\simeq z_2$ and $z_1,z_2\in A_1(k,n)$, then $z_2=\phi(z_1)$ with  $\phi\in  G'_0$.
	\end{prop}
	For a	nilpotent or semisimple element $z\in A_1$, we will call the element $z'\in A_1(k,n), k\ge n\ge 0,$ such that $z\simeq z'$   the Joseph normal form of $z$. In particular, the pair of integers $(k,n)$ with $k\ge n\ge 0$ is uniquely determined by $z$. 
	
If	$z$ is   nilpotent or weak semisimple, 
then its Joseph normal form  $z'$ is unique up to the action of the group $ G'_0$.

\begin{rem}
	As the field $K$ is assumed to be algebraically closed with characteristic 0, each strictly semisimple element of $A_1$ is conjugate to some $a\p\q+b,a\ne0$. Thus if $z$ is strictly semisimple then its Joseph normal form  $z'$ can be chosen as $ a\p\q+b,a\ne0$. Note that $\tau'(\p\q)=-(\p\q+1)$, and it is easy to verify that 	
	$a\p\q+b\simeq c\p\q+d$ iff $(c,d)=(a,b)$ or $(c,d)=(-a,b-a)$.

	Note that  each strictly nilpotent element $z\in A_1$  is conjugate to some element in $K[p]- K$, and  each strictly nilpotent element $z$  such that $C(z)=K[z]$ is conjugate to $\p$ \cite{d}.\\
\end{rem}

Recall that $S_1$ is  the 1st symplectic Poisson algebra. For $f\in S_1$, let \[ad_f:S_1\rt S_1, g\mapsto \{f,g\}.\]
 Note that  $ad_f$ is a $K$-derivation for $K[X,Y]$.  Let $Ev(ad_f)$ be the set of eigenvalues of $ad_f$ in $K$.  One can define $F(f), N(f), C(f), D(f)$ as in Section 6 of \cite{d}. Specifically, \[F(f)=\{g\in S_1|dim(K[ad_f]g)<\infty\},\]
 \[N(f)=\bigcup_{n>0} Ker\ (ad_f)^n,\ C(f)=Ker\ (ad_f),\]
 \[D(f)=\bigoplus_{\lambda\in Ev(ad_f)} Ker\ (ad_f- \lambda\cdot 1).\]
 
 One  has \[D(f)\cup N(f)\subseteq F(f), D(f)\cap N(f)=C(f).\] Note that $F(f), N(f), C(f), D(f)$ are all associative subalgebras of $K[X,Y]$. 
  As the GK-dimension of the $K$-algebra $K[X,Y]$ equals 2, the proof in Section 6 of \cite{d} applies to $S_1$. So analogous to the $A_1$ case, for all $f\in S_1-K$, either $F(f)=D(f)$ or $F(f)=N(f)$. See also Section 2 of \cite{b} for a proof of above result. It is clear that the center of the  Poisson algebra $S_1$ equals $K$.
      So just as in the $A_1$ case,  $S_1-K$ is a disjoint union of 5 nonempty sets $\Delta_i, i=1,\cdots,5,$ which is defined exactly as in  Section 6 of \cite{d}.  Elements in $\Delta_1\cup \Delta_2$ (resp. $\Delta_1$, $\Delta_2$) will be called nilpotent (resp. strictly nilpotent, weak nilpotent),  elements in $\Delta_3\cup \Delta_4$ (resp. $\Delta_3$, $\Delta_4$) will be called semisimple (resp. strictly  semisimple, weak  semisimple),  and elements in $\Delta_5$ will be called generic.
 
Recall that $\text{Aut}(S_1)$ is the automorphism group of the Poisson algebra $S_1$. One has $\tau\in \text{Aut}(S_1)$, and \[\text{Aut}(S_1)\cong \text{Aut}_n(K[X,Y]).\] 	For	$k,n\in \mathbb{Z}_{\geq 0}$, let $S_1(k,n)=\{f\in S_1|f = \sum a_{ij}X^i Y^j,a_{kn} \neq 0$, and if $ i> k$ or $j> n$ then $a_{ij} = 0\}$.
 
 The proof of Lemma 4.1 for $A_1$ in \cite{jo} also applies to $S_1$. So one has 
 
 \bco	Let	$f\in S_1.$ Assume $f$ is nilpotent or semisimple, then there exists $\phi \in \text{Aut}(S_1)$ such that $\phi(f)\in S_1(k,n)$, where $k,n\ge0, (k,n)\neq (0,0).$ 
 If $f$ is strictly semisimple, then $(k,n)=(1,1)$; otherwise $k\ne n$, and one can choose $k,n$ such that $k>n$.
	\eco
 For $f,g\in S_1$, we write $f\simeq g$ if for some $ \phi \in \text{Aut}(S_1), g= \phi(f)$.	Analogous to Proposition \ref{p1} one also has 
	
		\begin{prop}\lb{p2}
		Assume	$f\in S_1$  is nilpotent or semisimple, then there exists unique $(k,n)\in (\mathbb{Z}_{\geq 0})^2- \{ (0,0)\}, k\ge n\ge0$ such that there exists $\phi \in \text{Aut}(S_1)$ with $\phi(f)\in S_1(k,n)$. 	
		If $f$ is strictly semisimple, then $(k,n)=(1,1)$.		
		If $f$ is nilpotent or weak semisimple, then  $k>n$; if in addition $f\simeq f_1, f\simeq f_2$ and $f_1,f_2\in S_1(k,n)$, then $f_2=\phi(f_1)$ with  $\phi=(a_1X+c_1,a_2Y+c_2),\ a_1a_2=1$. 
		
	\end{prop}
		Analogously, if $f$ is strictly semisimple, then $f\simeq aXY+b,a\ne0$.
	One has  $\tau(XY)=-XY$ and it can be verified that 
	$aXY+b\simeq cXY+d$ iff $(c,d)=(a,b)$ or $(c,d)=(-a,b)$.		
	
		\section{Normal forms of elements and  Dixmier Conjecture}
	\setcounter{equation}{0}\setcounter{theorem}{0}
	Let	$\Gamma=\{z\in A_1|\exists w\in A_1, [z,w] = 1\}$, and for $z\in\Gamma$, $\Gamma(z) = \{w\in A_1| [z,w] = 1\}$. One knows that elements in 	$\Gamma$ are all nilpotent.
	
	
		For $z,w\in A_1$, let $<z,w>$ denote the $K$-subalgebra of $A_1$ generated by $z,w$.
		
	Note that $DC_1$ holds iff for $z,w\in A_1$ such that $[z,w] = 1$, $<z,w>=A_1$. Let us recall the following result before we prove an equivalent formulation of $DC_1$ in  Proposition \ref{p31}. 
	\bthm \lb{lm8}(Lemma 2.7 of \cite{d}, Proposition 3.2 of \cite{jo})\\
	Let $ z,w \in A_1- \{0\}$. Assume that $(\rho,\sigma)\in \mathbb{R}^2- \{(0,0)\}$ with $\rho+\sigma> 0$. Set $f=\mathbf{f}_{\rho,\sigma}(z)$ and $g=\mathbf{f}_{\rho,\sigma}(w)$. Then
	
	(1) $\mathbf{f}_{\rho,\sigma}(zw)=fg$;
	
	(2) If $\{f,g\}\ne 0$ then \be\lb{e9}  \mathbf{f}_{\rho,\sigma}([z,w])=-\{f,g\}\ee and $ \mathbf{v}_{\rho,\sigma}([z,w])=\mathbf{v}_{\rho,\sigma}(z)+\mathbf{v}_{\rho,\sigma}(w)-(\rho+\sigma);$
	
	(3) If $\{f,g\}= 0$, then $ \mathbf{v}_{\rho,\sigma}([z,w])<\mathbf{v}_{\rho,\sigma}(z)+\mathbf{v}_{\rho,\sigma}(w)-(\rho+\sigma)$.
	
	\ethm
	Note that there is a minus sign in (\ref{e9}) as our convention  $[\q,\p]=1$ is different to the convention in \cite{d}, which is $[\p,\q]=1$, while the identification (\ref{e00}) is the same as in \cite{d}.
	\begin{prop}\lb{p31}
		$DC_1$ holds $\Leftrightarrow $ For any $k,n$ with $k>n\geq 1$, $ A_1(k,n)\cap \Gamma=\o$.
	\end{prop}
	
	\bp
	$"\Rightarrow"$: Assume there exists $ k>n\geq 1,z\in A_1(k,n),z\in \Gamma$,
	$w\in \Gamma(z),[z,w] = 1$. Then \[\phi:A_1\rt A_1,\p\mapsto z,\q\mapsto -w\] is a homomorphism.
	If $<z,w>=A_1$, then $\phi \in \text{Aut}(A_1)$ and $\p\simeq z $, which contradicts to Corollary \ref{c3} as $\p\in A_1(1,0)$ and $z\in A_1(k,n), k>n\geq 1$. 
	So $<z,w>\subsetneqq A_1$, $DC_1$ is wrong.
	
	$"\Leftarrow"$: Assume $z\in\Gamma$. After applying some automorphism in  $\text{Aut}(A_1)$ one can assume that $z$ is in its Joseph normal form, i.e.
	$z\in A_1(k,n)\cap\Gamma$ for some  $k>n\geq 0$. By assumption,  if $k>n\geq 1$ then $A_1(k,n)\cap \Gamma=\o$, so now one must have $n = 0$, i.e. 
	$z\in A_1(k,0)\cap \Gamma$ with $k>0$.
	
	Choose $w\in \Gamma(z)$, such that $\mathbf{v}_{1,0}(w) = \inf\{\mathbf{v}_{1,0}(w') | w'\in \Gamma(z)\}$.
	Assume $\mathbf{f}_{1,0}(z) = a X^k =: f$, $a \neq 0$, and $\mathbf{f}_{1,0}(w) = X^ih(Y) =: g$, $i\geq 0$.
	
	If $\{f,g\} = 0$, then $g = b X^i, b\ne 0, i\ge0.$
	If $i>0$, then $ w' = w - b \p^i\in \Gamma(z),\mathbf{v}_{1,0}(w') < \mathbf{v}_{1,0}(w)$, which contradicts to the choice of $w$.
	So $i = 0$ thus $g= b$. Hence $ w=b$, which contradicts to $w\in \Gamma(z)$.
	Thus $\{f,g\}\ne0$, and by Lemma \ref{lm8},  $\{f,g\}=-\mathbf{f}_{1,0}([z,w]) $. 

So	 $\{f,g\}=-1$, i.e. $\{a X^k, X^ih(Y)\} = -1$. By direct computations one has $k = 1,i =0$. Now $z\in A_1(1,0)$ and $ z\simeq \p$, in this case one has $\Gamma(z)\subseteq K[\p]+K^{\times}\q $, thus  $ <z,w> = A_1$ for all $w\in \Gamma(z)$, and $DC_1$ holds.
	\ep

	\ble\label{lm3}(Proposition 3.1 of \cite{ht})
	Assume that $f\in K[X,Y]$ is $(\rho,\sigma)$-homogeneous with $|E(f)|>1$, where  $(\rho,\sigma)\ne(0,0)$ and $\rho+\sigma\ge 0$. If there exists $(\rho,\sigma)$-homogeneous $g\in K[X,Y]$ such that $\{f,g\}=1$, then $f$ equals one of the following:
	
	(1) $f_1=a X^n+b Y,  n\ge0,a,b\in K^\times$, and  $(\rho,\sigma)=(1,n)$; 
	
	(2)$f_2=\tau(f_1)=b X+a Y^n,  n\ge0,a,b\in K^\times$, and $(\rho,\sigma)=(n,1)$.

	\ele

	\ble\label{lm4}
	Assume   $\rho+\sigma>0$. Assume  $z\in A_1$,  $\mathbf{f}_{\rho,\sigma}(z)=f$, and $|E(f)|>1$. If
	(1) $f$ is not a proper power in $K[X,Y]$ (i.e. $f$ cannot be written as $f= g^n$ with $ g\in K[X,Y]$ and $n>1$), and	
	(2)  $f$ is not of the type of polynomials in Lemma \ref{lm3}, then $z\notin \Gamma$.
	\ele
	\bp
	Assume that $z\in \Gamma$ and $w\in \Gamma(z)$.
	
	If $\mathbf{v}_{\rho,\sigma}(z)\le0$, then one of  $\rho,\sigma<0$. Assume $\rho<0$ (the case $\sigma<0$ is similar), then $z \in \sum_{(i,j):i\geq j\geq 0}K\p^i\q^j$.  As $[z,w]=1$, by Theorem 3.7 of \cite{ht}, $z=a \p+b,a\neq 0$ and $|E(f)|=1$, contradicts to $|E(f)|>1$. So one must have $\mathbf{v}_{\rho,\sigma}(z)>0$.
	
	Assume $\mathbf{f}_{\rho,\sigma}(w) = g$.
	As $[z,w] = 1$, if $\{f,g\} \neq 0$ then $\{f,g\} =-\mathbf{f}_{\rho,\sigma}([z,w])=-1$ , 
	which contradicts to (2). 
	
	So $\{f,g\} = 0$. Then by  (1), there exists $ m_0 \geq 0, g\sim f^{m_0}$. If $m_0>0$, then $g= a_0f^{m_0}$
	for some $ a_0 \in K^{\times}$. Set $w_1 = w - a_0 z^{m_0}$. Then  $[z,w_1] = 1$  and
	$\mathbf{v}_{\rho,\sigma}(w_1) < \mathbf{v}_{\rho,\sigma}(w)$. Let  $\mathbf{f}_{\rho,\sigma}(w_1) = g_1$. By the same reason as above, one has 
	$\{f,g_1\} = 0$, and there exists $ m_1 \geq 0, g\sim f^{m_1}$. If $m_1>0$ then we continue this process until we find some  $w' = w - \sum_{i = 0}^{s}a_i z^{m_i} \in \Gamma(z)$, $\mathbf{f}_{\rho,\sigma}(w') = h$, and $ h\sim f^0=1$.
	
	If $\rho,\sigma\ge0$, then $h$, thus $w'$, is in $K^\times$, which contradicts to  $ w\in \Gamma(z)$. If one of $\rho,\sigma<0$, for example $\rho<0$ (the case $\sigma<0$ is similar), then $w \in \sum_{(i,j):i\geq j\geq 0}K\p^i\q^j$. Then as $[z,w]=1$, by Theorem 3.7 of \cite{ht}, $w=a \p+b,\ a,b\ne0;$  $z=a^{-1} \q+ h(\p), h\in K[X]$. Then $ f(X,Y)=a^{-1} Y$ and  $|E(f)|=1$, contradicts to $|E(f)|>1$.
	
	\ep

	\begin{prop}\label{prop4} Let $F$ be a field of characteristic 0. 
		Assume	$f(X) = (X-a_1)^{l_1}\cdots(X-a_k)^{l_k}\in F[X]$, $k\geq 2$, $l_i\geq 1,a_i\in F\ for\ i=1,\cdots,k$, and $a_i\ne a_j$ for $i\ne j$. One has the partial fraction decomposition
		$$\dfrac{1}{f(X)}=\sum_{i=1}^{k}(\dfrac{c_{i,1}}{X-a_i}+\cdots+\dfrac{c_{i,l_i}}{(X-a_i)^{l_i}}), c_{i,j}\in F.$$
		Then there exists some $i$ such that $c_{i,1}\neq 0$. 		
	\end{prop}
	The following proof is given by Xiaoguang Wang.
	\begin{proof}
		If $l_j=1$ for some $j$, then the conclusion follows trivially. So we assume that  $l_j>1$ for all $j$.
		
		Assume $c_{i,1}=0$ for $i=1,\cdots,k$, then 
		\begin{align}
			\dfrac{1}{f(X)}&=\sum_{i=1}^{k}(\dfrac{c_{i,2}}{(X-a_i)^2}+\cdots+\dfrac{c_{i,l_i}}{(X-a_i)^{l_i}})\notag\\
			&=\sum_{i=1}^{k}(\dfrac{-c_{i,2}}{X-a_i}+\cdots+\dfrac{-c_{i,l_i}}{(l_i-1)(X-a_i)^{l_i-1}})'\notag\\
			&=(\dfrac{g(X)}{h(X)})', \label{eq3}
		\end{align}
		where $g,h\in F[X]$, and $h(X)=\Pi_{i=1}^k (X-a_i)^{l_i-1}$. It is clear that		
		$\deg g<\deg h$ and \be
		\deg h =\sum_{i=1}^{k}(l_i-1)=\deg f-k
		\label{eq1}.\ee By (\ref{eq3}),
		\be
		f(g'h-gh')=h^2.
		\label{eq2}\ee
		As $\deg g< \deg h$, the degree of LHS of (\ref{eq2}) equals  $\deg f+\deg g+\deg h -1$. The degree of RHS is $2\deg h$. Thus 
		$$\deg f+\deg g-\deg h -1=0.$$
		Substitute  (\ref{eq1}) into this equation, one has
		$$\deg g+k-1=0,$$
		which contradicts to $ k\geq 2$.
	\end{proof}
		\ble\lb{lm9}
	Assume $f\in K(X)$ and $ a\in K$ is a pole of $f$ of order 1. Then there exists no $g\in K(X)$ such that $g'=f$.
	\ele
	\bp
	Assume $g'=f$ for some $g\in K(X)$. If $a$ is a pole of $g$, say, of order $k\ge1$, then 
	$a$ is a pole of $f$ of order $k+1>1$, which contradicts to the assumption. If $a$ is not a pole of $g$, then  $a$ is not a pole of $f$ either,  which also  contradicts to the assumption.	
	\ep
	
	For $0\ne h\in K[X]$, let $r(h)$ denote the number of distinct roots of $h$.
	\ble\label{lm6}
	Assume
	$ k> n \geq 2$ and $z\in  A_1(k,n)$. 
	
	(1) If $\mathbf{f}_{0,1}(z)$ is not a proper power in $K[X,Y]$, then $z\notin \Gamma$.
	
	(2) If $\mathbf{f}_{0,1}(z)$ is a proper power and equals $(h_0(X)Y)^n,n\geq 2$, $h_0\in K[X]$ and $r(h_0)>1$, then  $z\notin \Gamma$.
	\ele
	\bp
	Assume $z\in \Gamma$. We will show that it leads to a contradiction under the hypothesis (1) or (2).
	Let $ w\in \Gamma(z)$, such that
	\begin{equation}\label{l2}
		\mathbf{v}_{0,1}(w) = \inf\{\mathbf{v}_{0,1}(\tilde{w}) | \tilde{w} \in \Gamma(z)\},
	\end{equation} and let $\mathbf{v}_{0,1}(w) = l$. One has $l\ge0$.
	
	Assume $\mathbf{f}_{0,1}(z) = f(X,Y) = h(X)Y^n, deg(h) = k$, and $\mathbf{f}_{0,1}(w) = g(X,Y)=h^*(X)Y^l, h^*\ne0$.
	As $[z,w] = 1$, if $\{f,g\} \neq 0$ then $\{f,g\} =-\mathbf{f}_{0,1}([z,w])= -1$ and $0 = \mathbf{v}_{0,1}([z,w]) = n + l - 1$,
	which contradicts to $n\ge2$. So $\{f,g\} = 0$.
	
	(1) Assume $\mathbf{f}_{0,1}(z)$ is not a proper power in $K[X,Y]$.
	
	Then by $\{f,g\} = 0$, there exists $ m \geq 0, g\sim f^m$ and $g=a f^m$ for some $ a \in K^{\times}$. If $m>0$, then
	set $w' = w - a z^m$. One has $[z,w'] = 1$ and
	$\mathbf{v}_{0,1}(w') < \mathbf{v}_{0,1}(w)$, which contradicts to (\ref{l2}).  If $m=0$ then $g$, thus $w$, is some nonzero constant in $K$, which contradicts to  $ w\in \Gamma(z)$. Thus (1) of this lemma is proved.
	
	(2) Assume $\mathbf{f}_{0,1}(z)$ is a proper power and equals $(h_0(X)Y)^n$, $ h_0\in K[X]$ and $ r(h_0)>1$. 
	Recall that \[\widehat{D}_{0,1} = \{\sum_{i\leq r}a_i(\p)\q^i |r\in \mathbb Z, a_i \in K(X), a_r\ne0\}\cup \{0\} \] is the completion of the quotient division algebra $Q(A_1)$ of $A_1$  with respect to $\mathbf{v}_{0,1}$; see \cite{hcp}\cite{g}. One has $z,w \in A_1 \subset \widehat{D}_{0,1}$. As
	$\mathbf{f}_{0,1}(z) = (h_0(X)Y)^n$, by Proposition 2.7 of  \cite{g}, there exists some unique $z_0\in \widehat{D}_{0,1}$,
$ \mathbf{f}_{0,1}(z_0) = h_0(X)Y:=f_0(X,Y)$ and 	$z = z_0^n$.
	
	As $\mathbf{v}_{0,1}(z)+\mathbf{v}_{0,1}(w)-1>0$,
	$  \{f_0^n,g\}=\{f,g\} =0$. Since $f_0=h_0(X)Y$ is not a proper power in $K(X,Y)$, $g\sim f_0^{l_0}$ for some  $l_0>0$. Then $g=b_0 f_0^{l_0}$ for some
	$ b_0 \in K^{\times}$. Set $w_1 = w - b_0z_0^{l_0}$, then $[z,w_1] = 1$ and $\mathbf{v}_{0,1}(w_1) < \mathbf{v}_{0,1}(w)$.
	
	If $\mathbf{v}_{0,1}(w_1) > 1- \mathbf{v}_{0,1}(z)$, then let $l_1 = \mathbf{v}_{0,1}(w_1), g_1 = \mathbf{f}_{0,1}(w_1)$. Then one still has 
	$\{f,g_1\} = \{f_0^n , g_1\}=0$, $ g_1 \sim f_0^{l_1}$ and  $ g_1=b_1 f_0^{l_1}$ for some $b_1 \in K^{\times}$.
	Set $w_2 = w_1 - b_1z_0^{l_1}$, then $[z,w_2] = 1$ and $\mathbf{v}_{0,1}(w_2) < \mathbf{v}_{0,1}(w_1)$.
	
	Continue in this way, one can find $w' = w - \sum_{i = 0}^{s}b_i z_0^{l_i}$,
	such that \[ [z,w'] = 1, \mathbf{v}_{0,1}(w') = 1- \mathbf{v}_{0,1}(z) = 1-n.\]
	Let $\mathbf{f}_{0,1}(w') = \bar{g}(X,Y)= r(X)Y^{1-n}\in K(X,Y)$ with $ r\in K(X)$. Then $\{f,\bar{g}\} = - \mathbf{f}_{0,1}([z,w']) =-1$, by Theorem \ref{lm8} (which is also valid in $\widehat{D}_{0,1}$). So
	\[-1 = \{f,\bar{g}\} = \{h(X)Y^n, r(X)Y^{1-n}\} =  (1-n)h'(X)r(X)-nh(X)r'(X).\]
Dividing $-nh(X)$ on both sides, one has \[ r' + (\frac{n-1}{n}\frac{h'}{h})r = \frac{n^{-1}}{h}.\]
	As $h = h_0^n$, one has in $K(X)$,
	\begin{align}\nonumber
		&\ \ \  r' + (\frac{n-1}{n}\frac{h'}{h})r = \frac{n^{-1}}{h},\\\nonumber
		&\Leftrightarrow r' + (n-1)\frac{h_0'}{h_0}\cdot r = \frac{n^{-1}}{h},\\\nonumber
		&\Leftrightarrow h_0^{n-1}[r' + (n-1)\frac{h_0'}{h_0}\cdot r] = \frac{n^{-1}}{h_0},\\\nonumber
		&\Leftrightarrow (n\cdot h_0^{n-1} r)' = \frac{1}{h_0}.\\    
	\end{align}  

	As $K$ is algebraically closed, one can assume $h_0(X) = b\prod_{i = 1}^{t}(X - \alpha_i)^{m_i}, b\ne0,  m_i \geq 1, \alpha_i \neq \alpha_j $ for $ i \neq j$.
	Then $r(h_0)=t$. By assumption $t\ge2$.
	
	Then one has the partial fraction decomposition 
	\[\frac{1}{h_0(X)} = \sum_{j = 1}^{t}\sum_{u = 1}^{m_j}\frac{b_{j,u}}{(X-\alpha_j)^u},\ b_{j,u}\in K.\] By Proposition \ref{prop4}, there exist some $ j$'s with $ b_{j,1} \neq 0$.
	Then  $\frac{1}{h_0(X)}=F_1(X) + F_2(X)$ where \[0 \neq F_1(X) := \sum_{j:b_{j,1}\neq 0}\frac{b_{j,1}}{X-\alpha_j},\  F_2(X)=\frac{1}{h_0(X)}-F_1(X).\]
	Obviously there exists $ \tilde{F_2} \in K(X)$  such that $\tilde{F_2}' = F_2$, then
	$ ( n h_0^{n-1} r-\tilde{F_2})'=F_1 $. But as $F_1$ has poles of order 1, this contradicts to Lemma \ref{lm9}.
	Thus (2) of this lemma is proved.
	
	\ep
	
	\ble\label{lem1} If
	$k> n \geq 1$ and $ gcd(k,n) = 1$, then
	$\Gamma\cap A_1(k,n)=\o$. In particular, 	$\Gamma\cap A_1(k,1)=\o$ for $k>1$.
	\ele
	\bp
	Assume
	$k> n \geq 1, gcd(k,n) = 1$, and
	$z\in\Gamma\cap A_1(k,n)$.  Let $w\in \Gamma(z)$ satifying  
	$\mathbf{v}_{1,1}(w) = \inf\{\mathbf{v}_{1,1}(\tilde{w}) | \tilde{w} \in \Gamma(z)\}$. 
	Then $\mathbf{f}_{1,1}(z) \sim X^kY^n$, and we will assume that 
	$\mathbf{f}_{1,1}(z)=X^kY^n$ without loss of generality.  Let $\mathbf{f}_{1,1}(w) = g(X,Y)$. 
	
	If $\{X^kY^n, g(X,Y)\} = 0$, then as $gcd(k,n) = 1$, there exists $ l \geq 0, g(X,Y)\sim (X^kY^n)^l$ and $g(X,Y)=a (X^kY^n)^l$ for some $ a \in K^{\times}$. If $l>0$ then set $ w' = w - a z^l$, one has $[z, w'] = 1$ and $\mathbf{v}_{1,1}(w') < \mathbf{v}_{1,1}(w)$,
	contradicting to $\mathbf{v}_{1,1}(w)$ being minimal. If $l=0$ then $w\in K^\times$, which contradicts to $w\in \Gamma(z)$.
	So $\{X^kY^n, g(X,Y)\} \neq 0$. Then 
	\begin{equation}
		\mathbf{v}_{1,1}(z) + \mathbf{v}_{1,1}(w) -(1+1) = 0\label{e1}
	\end{equation}
	and $\{X^kY^n, g(X,Y)\} =-\mathbf{f}_{1,1}([z,w])= -1$.
	As $\mathbf{v}_{1,1}(w) \geq 0$, by (\ref{e1}), $k+n \leq 2$, which contradicts to $k>n\geq 1$.
	\ep
	
	\bthm \lb{main}Assume
	$ k> n \geq 1$. If $k$ or $n$ is prime, then $ \Gamma\cap A_1(k,n)=\o$.
	\ethm
	\bp
	If $k$ is prime, then as $k> n \geq 1$, one has $gcd(k,n) = 1$.
	By Lemma \ref{lem1}, $ \Gamma\cap A_1(k,n)=\o$. In the following we assume $n$ is prime and write $n = p$. 
	
	Assume  $ \Gamma\cap A_1(k,p)\ne\o$. (We will show that it will lead to a contradiction.) By  Lemma \ref{lem1}, one has $p|k$ (and $p<k$).
	
	Assume $z\in \Gamma\cap A_1(k,p)$.
	
	If $\mathbf{f}_{0,1}(z)$ is not a proper power in $K[X,Y]$, then by Lemma \ref{lm6} (1), $z\notin\Gamma$, contradicts to the assumption that $z\in\Gamma$.
	
	So $\mathbf{f}_{0,1}(z)$ is a proper power in $K[X,Y]$. Assume $\mathbf{f}_{0,1}(z) = h(X)Y^p$.
	As $p$ is prime, there exists some $ h_0\in K[X], h = h_0^p$ and $\mathbf{f}_{0,1}(z) = (h_0(X)Y)^p$. 
	
	Assume $h_0$ has $t$ different roots. If $t\geq 2$, then by Lemma \ref{lm6} (2), $z\notin\Gamma$, contradicts to the assumption that $z\in\Gamma$. Hence $t = 1$.
	
	Now one has $h_0(X) \sim (X-\alpha)^{k_0}, k_0 \geq 1,\alpha\in K$. So $\mathbf{f}_{0,1}(z) = a(X-\alpha)^kY^p$, $k = k_0p$, and $a \in K^{\times}$.
	Let $\phi \in \text{Aut}(A_1)$ such that  $\phi$ maps $\p\mapsto \p+\alpha, \q\mapsto \q$. Replace $z$ by $a^{-1}\phi (z)$, one still has  $z\in \Gamma\cap A_1(k,p)$. Now $\mathbf{f}_{0,1}(z) = X^kY^p$.
	
	Assume there exists $ \sigma\in \mathbb R$, such that $\sigma > k_0$, $X^kY^p\in \mathbf{f}_{-1,\sigma}(z)$, and $|E(\mathbf{f}_{-1,\sigma}(z))| \geq 2$.
	
	If $\mathbf{f}_{-1,\sigma}(z)$ is a proper power in $K[X,Y]$, as $p$ is prime, $\mathbf{f}_{-1,\sigma}(z) = f_0(X,Y)^p$ implies that $ f_0(X,Y) = b X^{k_0}Y$ is a monomial. (See Figure 1.)
	So $\mathbf{f}_{-1,\sigma}(z)= f_0(X,Y)^p=(b X^{k_0}Y)^p=b^p X^kY^p$ is a monomial,  which contradicts to  $|E(\mathbf{f}_{-1,\sigma}(z))| \geq 2$.
	Thus $\mathbf{f}_{-1,\sigma}(z)$ is not a proper power in $K[X,Y]$. By Lemma \ref{lm4}, $z\notin\Gamma$, which is a contradiction.
	\begin{figure} [htp]
		\centering
		\begin{tikzpicture}
			\draw [help lines, gray] (0,0) grid (8,6);
			\draw [->, thick](0,0) -- (9,0);
			\draw [->, thick](0,0) -- (0,7);
			\draw [->, thick](2,3) -- (1,6);
			\draw [thick](0,1) -- (6,3);
			\draw [fill] (6,3) circle [radius = 1.5pt];
			\node[right] at (6,3) {$X^kY^p$};
			\draw [fill] (3,2) circle [radius = 1.5pt];
			\node[above] at (2.5,4) {$(-1,\sigma)$};
			\draw [fill] (2,1) circle [radius = 1.5pt];
			\node[right] at (2,1) {$X^{k_0}Y$};
			\draw [very thick] (6,0) -- (6,3);
			\draw [dashed] (0,1/3) -- (2,1);
			\draw [dashed] (2,0) -- (2,1);
		\end{tikzpicture}\\
		{Figure 1}
	\end{figure}

	Thus  there exists no $ \sigma\in \mathbb R$ such that $\sigma > k_0$, $X^kY^p\in \mathbf{f}_{-1,\sigma}(z)$, and $|E(\mathbf{f}_{-1,\sigma}(z))| \geq 2$. So 
	$z \in \sum_{(i,j):i\geq j\geq 0}K\p^i\q^j$. Then as $z\in\Gamma$, by Theorem 3.7 of \cite{ht}, $z=a\p+b\in A_1(1,0), a\ne0$, contradicts to $ z\in A_1(k,p)$.
	
	\ep

	\begin{rem} By \cite{ggv},  
		if	$z\in A_1(k,n), k> n \geq 1, k+n\le15$, then $[z,w] = 1$ has no solution for $w \in A_1$. So $ \Gamma\cap S_1(k,n)=\o,  k> n \geq 1, k+n\le15$.	
	\end{rem}
	
	Now we state the parallel results for  the Poisson algebra $S_1$ and 
	$JC_2$.
		Let	$$\Gamma'=\{f\in S_1|\exists g\in S_1, \{f,g\} = 1\},$$ and for $f\in\Gamma'$, $\Gamma'(f) = \{g\in S_1| \{f,g\} = 1\}$. Elements in $\Gamma'$ are all nilpotent.
		\ble\lb{}
		$JC_2$ holds $\Leftrightarrow$ If $\{f,g\}=1$ then $(f,g)\in \text{Aut}(K[X,Y])$.
	\ele  This is obvious.
		\begin{prop}\lb{}
		$JC_2$ holds $\Leftrightarrow $ If $k>n\ge1$ then $S_1(k,n)\cap \Gamma'=\o$.
	\end{prop}
	\bp
		$"\Rightarrow"$: Assume there exist $k,n$ with $ k>n\geq 1,f\in S_1(k,n)\cap \Gamma'$, $ \{f,g\}=1$ for some $g\in  S_1$. As $JC_2$ is assumed to be true,  $(f,g)\in \text{Aut}(K[X,Y])$. As $f$ is of rectangular type, by Proposition \ref{th1}, $f=a_1X+a_0\ or\ a_1Y+a_0, a_1\ne 0$. Then $f\in S_1(1,0)\ or\  S_1(0,1)$, which contradicts to $ f\in S_1(k,n),k>n\geq 1$.
	
	$"\Leftarrow"$: Assume $\{f,g\}=1$. Then $f$ is nilpotent in $S_1$, and one can assume that $f$ is of rectangular type $(k,n)$ with $k>n\ge0$.
	 By assumption, $S_1(k,n)\cap \Gamma'=\o$ for  $k>n\ge1$, so $f\in S_1(k,0)$ for some $k>0$. Thus $f$ is a polynomial in $X$. By $\{f,g\}=1$ one has  $f=a_1X+a_0, a_1\ne0$, and $g=a_1^{-1}Y+h(X), h\in K[X]$, thus  $(f,g)\in \text{Aut}(K[X,Y])$.
	
	\ep
	
	Lemma \ref{lem1} and Theorem \ref{main} can also be generalized to the  $S_1$ case, so one has 
	
		\bco \lb{}Assume
	$ k> n \geq 1$. If $gcd(k,n)=1$, or one of $k$, $n$ is prime, then $ \Gamma'\cap S_1(k,n)=\o$.
	\eco

	
\end{document}